\documentclass[11pt]{article}
	\usepackage{amsfonts}
	\newcommand{\blind}{0}
	
    \makeatletter
    \renewcommand\section{\@startsection {section}{1}{\z@}%
                                       {-3.5ex \@plus -1ex \@minus -.2ex}%
                                       {2.3ex \@plus.2ex}%
                                       {\normalfont\fontfamily{phv}\fontsize{16}{19}\bfseries}}
    \renewcommand\subsection{\@startsection{subsection}{2}{\z@}%
                                         {-3.25ex\@plus -1ex \@minus -.2ex}%
                                         {1.5ex \@plus .2ex}%
                                         {\normalfont\fontfamily{phv}\fontsize{14}{17}\bfseries}}
    \renewcommand\subsubsection{\@startsection{subsubsection}{3}{\z@}%
                                        {-3.25ex\@plus -1ex \@minus -.2ex}%
                                         {1.5ex \@plus .2ex}%
                                         {\normalfont\normalsize\fontfamily{phv}\fontsize{14}{17}\selectfont}}
    \makeatother
	
	\usepackage{hyperref}
        \usepackage{amsmath}
	\usepackage{graphicx}
	\usepackage{enumerate}
	\usepackage{xcolor}
        \usepackage{natbib}
	\usepackage{url} 
        \newtheorem{definition}{Definition}[section]
        
\begin{document}
		
		\def\spacingset#1{\renewcommand{\baselinestretch}%
			{#1}\small\normalsize} \spacingset{1}
		
		\if0\blind
		{
			\title{\bf \textcolor{black}{Survey on} Lagrangian Relaxation for MILP: Importance, Challenges, \textcolor{black}{Historical Review,} Recent Advancements, and Opportunities}
\author{Mikhail A. Bragin\\
Department of Electrical and Computer Engineering, \\
University of Connecticut, Storrs, USA\\
\texttt{mikhail.bragin@uconn.edu
}}

			\date{}
			\maketitle
		} \fi
		
		\if1\blind
		{

            \title{\bf \emph{IISE Transactions} \LaTeX \ Template}
			\author{Author information is purposely removed for double-blind review}
			
\bigskip
			\bigskip
			\bigskip
			\begin{center}
				{\LARGE\bf \emph{IISE Transactions} \LaTeX \ Template}
			\end{center}
			\medskip
		} \fi
		\bigskip
		
	\begin{abstract}
\noindent Operations in areas of importance to society are frequently modeled as Mixed-Integer Linear Programming (MILP) problems. While MILP problems suffer from combinatorial complexity, Lagrangian Relaxation has been a beacon of hope to resolve the associated difficulties through decomposition. Due to the non-smooth nature of Lagrangian dual functions, the coordination aspect of the method has posed serious challenges. This paper presents several significant historical milestones (beginning with Polyak's pioneering work in 1967) toward improving Lagrangian Relaxation coordination through improved optimization of non-smooth functionals.
Finally, this paper presents the most recent developments in Lagrangian Relaxation for fast resolution of MILP problems. The paper also briefly discusses the opportunities that Lagrangian Relaxation can provide at this point in time.
	\end{abstract}
			
	\noindent%

	\spacingset{1.5} 

\section{Introduction} \label{sec1}

The aim of this paper is to review Lagrangian-Relaxation-based methods for Mixed-Integer Linear Programming (MILP) problems. Because of the integer variables, Lagrangian Relaxation leads to non-smooth optimization in the dual space. Accordingly, key non-smooth optimization methods will also be reviewed. 

This paper focuses on the efficient resolution of \textit{separable} Mixed-Integer Linear Programs (MILPs), which are formally defined as follows:

\begin{flalign}
& \min_{(x,y) := \left\{x_i,y_i\right\}_{i=1}^I} \Bigg\{\sum_{i=1}^I \left((c_i^x)^T \cdot x_i + (c_i^y)^T \cdot y_i\right)  \Bigg\},  \label{eq1}
\end{flalign}
whereby $I$ subsystems are coupled through the following constraints  
\begin{flalign}
& s.t. \;\; \sum_{i=1}^I A_i^x \cdot x_i + \sum_{i=1}^I A_i^y \cdot y_i - b = 0, \;\; \left\{x_i,y_i\right\} \in \mathcal{F}_i, i = 1, \dots, I.  \label{eq2}
\end{flalign}

\noindent The \textit{primal problem} \eqref{eq1}-\eqref{eq2} is assumed to be feasible and the feasible region $\mathcal{F} \equiv \prod_{i=1}^I \mathcal{F}_i$ with $ \mathcal{F}_i \subset \mathbb{Z}^{n_i^x} \times \mathbb{R}^{n_i^y}$ is assumed to be bounded and finite.


\subsection{Importance and Difficulties of MILP Problems}
MILP has multiple applications in problems of importance to  society: 
ambulance relocation \citep{lee2022lagrangian}, 
cost-sharing for ride-sharing \citep{hu2021cost}, drop box location \citep{schmidt2022locating}, efficient failure detection in large-scale distributed systems \citep{er2022data}, home healthcare routing \citep{dastgoshade2020lagrangian}, home service routing and appointment scheduling \citep{tsang2022stochastic}, 
job-shop scheduling \citep{liu2021novel}, facility location \citep{BASCIFTCI2021548}, flow-shop scheduling \citep{hong2019admission, balogh2022milp, oztop2022metaheuristics}, freight transportation \citep{archetti2021optimization}, location and inventory prepositioning of disaster relief supplies \citep{shehadeh2022stochastic}, machine scheduling with sequence-dependent setup times \citep{yalaoui2021identical}, maritime inventory routing \citep{gasse2022machine}, multi-agent path finding with conflict-based search \citep{huang2021learning}, multi-depot electric bus scheduling \citep{gkiotsalitis2022exact}, multi-echelon/multi-facility green reverse logistics network design \citep{reddy2022multi}, optimal physician staffing \citep{prabhu2021physician}, optimal search path with visibility \citep{morin2023ant}, oral cholera vaccine distribution \citep{smalley2015optimized}, outpatient colonoscopy scheduling \citep{shehadeh2020distributionally}, pharmaceutical distribution \citep{zhu2018design}, plant factory crop scheduling \citep{huang2020plant}, post-disaster blood supply \citep{hamdan2020robust,kamyabniya2021robust}, real assembly line balancing with human-robot collaboration \citep{nourmohammadi2022balancing}, reducing vulnerability to human trafficking \citep{kaya2022improving}, restoration planning and crew routing \citep{MORSHEDLOU2021107626}, ridepooling \citep{gaul2022event}, scheduling of unconventional oil field development \citep{soni2021mixed}, security-constrained optimal power flow \citep{velloso2021exact}, semiconductor manufacturing \citep{chang2017stochastic}, surgery scheduling \citep{kayvanfar2021new}, unit commitment \citep{kim2018temporal, chen2019distributed, li2019multi, chen2020high, li2020robust, van2021decomposition}, vehicle sharing and task allocation \citep{arias2022vehicle}, workload apportionment \citep{gasse2022machine}, and many others.  
However, MILP problems belong, \textcolor{black}{in general,} to the class of NP-hard problems because of the presence of integer variables $x$.  MILP problems of practical sizes are generally difficult to solve due to their combinatorial complexity. As the problem size increases, the computational effort required to obtain an optimal solution increases superlinearly, e.g., exponentially; for a number of practical problems, the computational effort may be significant to obtain even a feasible solution. \textcolor{black}{Almost all optimization algorithms thus have a superlinear running time given the NP-hard nature of these problems, which implies that there is no known polynomial-time algorithm to solve them optimally unless P=NP.} 
Moreover, many problems of importance typically require short solving times (ranging from 20 minutes to a few seconds, depending on the \textcolor{black}{application}), as well as high-quality solutions.

With decreasing problem size, the NP-hardness has the property of reducing complexity \textcolor{black}{exponentially}. The \textit{dual} decomposition and coordination Lagrangian Relaxation method is promising to exploit this reduction of complexity; the method essentially ``reverses'' combinatorial complexity upon decomposition, thereby drastically reducing the effort required to solve subproblems (each subproblem $i$ corresponds to a subsystem $i$). Lagrangian Relaxation is also deeply rooted in economic theory, whereby the solutions obtained are rested upon the economic principle of ``supply and demand.'' When the ``demand'' exceeds the ``supply,'' Lagrangian multipliers (which can be viewed as ``shadow prices'') increase (and vice versa) thereby discouraging subsystems from making less ``economically viable'' decisions. Notwithstanding the advantage of the decomposition aspect, the ``price-based'' coordination of the method (to appropriately coordinate the subproblems), however, has been the subject of intensive research for many decades because of the fundamental difficulties of the underlying non-smooth optimization of the associated dual functions caused by the presence of integer variables in the primal space. Accordingly, key non-smooth optimization methods will also be reviewed. 

The purpose of this paper is to present a brief overview of the key milestones in the development of the Lagrangian Relaxation method for MILP problems as well as in the optimization of convex non-smooth functions. The rest of the paper is structured as follows: 
\begin{enumerate}
   \item At the beginning of Section \ref{sec2}, the Lagrangian dual problem 
   is presented and the difficulties of Lagrangian Relaxation on a pathway to solving MILP problems are explained. In subsequent subsections, the difficulties are resolved one by one; 
   \item In subsection \ref{sec21}, early research on non-smooth optimization \citep{polyak1967general, polyak1969minimization} is presented to lay the foundation for further developments; specifically, the Polyak formula \citep{polyak1969minimization} depending on the \textit{optimal dual value} $q(\lambda^*)$ to ensure \textit{geometric} (also referred to as \textit{linear} convergence rate is presented; 
   
      \item \textcolor{black}{In subsection \ref{sec22a}, fundamental research as well as applications of Lagrangian duality that emerged in 1970's to solve discrete optimization problems is discussed;}
   \item In subsection \ref{sec22}, the \textit{subgradient-level} method \citep{goffin1999convergence} is presented to ensure convergence without the need to know $q(\lambda^*)$;

    \item In subsection \ref{sec23}, the fundamental difficulties associated with subgradient methods (high computational effort and zigzagging of multipliers) are explained; 
    
    \item In subsections \ref{sec24} and \ref{sec25}, two separate research thrusts: \textit{surrogate} \citep{kaskavelis1998efficient, zhao1999surrogate} and \textit{incremental} \citep{nedic2001incremental}) to reduce computational effort as well to alleviate zigzagging of multipliers are reviewed; the former thrust still requires $q(\lambda^*)$ for convergence; the latter thrust avoids the need to know $q(\lambda^*)$ following the ``subgradient-level'' ideas presented in subsection \ref{sec22}; \item In subsection \ref{sec26}, Surrogate Lagrangian Relaxation (SLR) \citep{bragin2015convergence} that proved convergence without $q(\lambda^*)$ by exploiting ``contraction mapping'' while inheriting convergence properties of the \textit{surrogate} method of subsection \ref{sec24} is reviewed; 
    \item In subsection \ref{sec27}, further methodological advancements for SLR are presented; to accelerate the reduction of constraint violations while enabling the use of MILP solvers, ``absolute-value'' penalties have been introduced \citep{bragin2018scalable}; to efficiently coordinate distributed entities while avoiding the synchronization overhead, computationally distributed version of SLR has been developed to efficiently coordinate distributed subsystems in an asynchronous way \citep{bragin2020distributed};
    \item In subsection \ref{sec28}, Surrogate ``Level-Based'' Lagrangian Relaxation \citep{bragin2022surrogate} is reviewed; this first-of-the-kind method exploits the linear-rate convergence potential intrinsic to the Polyak's formula presented in subsection \ref{sec21} but without the knowledge $q(\lambda^*)$ and without heuristic adjustments of its estimates presented in subsection \ref{sec22}. Rather, an estimate of $q(\lambda^*)$ (the ``level value'') has been innovatively determined purely through a simple constraint satisfaction problem; the surrogate concept \citep{zhao1999surrogate, bragin2015convergence} ensures low computational requirements as well as the alleviated zigzagging; accelerated reduction of constraint violations is achieved through ``absolute-value'' penalties \citep{bragin2018scalable} enabling the use of MILP solvers;
    \item In Section \ref{Conclusion}, a brief conclusion is provided with future directions delineated. 
\end{enumerate}

\section{Lagrangian Duality for Discrete Programs and Non-Smooth Optimization} \label{sec2}
The Lagrangian \textit{dual problem} that corresponds to the original  MILP problem \eqref{eq1}-\eqref{eq2} is the following \textit{non-smooth} optimization problem:
\begin{flalign}
& \max_{\lambda} \{q(\lambda): \lambda \in \Omega \subset \mathbb{R}^m\},  \label{eq3}
\end{flalign}
where the convex \textit{dual function} is defined as follows:
\begin{flalign}
& q(\lambda) = \min_{(x,y)} \big\{L(x,y,\lambda), \left\{x_i,y_i\right\} \in \mathcal{F}_i, i = 1,\dots,I \big\}.  \label{eq4}
\end{flalign}
Here $L(x,y,\lambda) \equiv \sum_{i=1}^I \big((c_i^x)^T \cdot x_i + (c_i^y)^T \cdot y_i\big) + \lambda^T \cdot \big(\sum_{i=1}^I A_i^x \cdot x_i + \sum_{i=1}^I A_i^y \cdot y_i - b\big)$ is the \textit{Lagrangian function} obtained by relaxing coupling constraints \eqref{eq2} by using Lagrangian multipliers $\lambda$. The minimization of $L(x,y,\lambda)$ within \eqref{eq4} with respect to $\{x,y\}$ is referred to as the \textit{relaxed problem}, which is separable into subproblems due to the additivity of $L(x,y,\lambda)$. This feature will be exploited starting from subsection \ref{sec24}. 

Even though the original \textit{primal problem} \eqref{eq1} is non-convex, $q(\lambda)$ is always continuous and concave with the feasible set $\Omega$ being always convex. Due to integer variables $x$ in the primal space, $q(\lambda)$ is non-smooth with facets (each representing a particular solution to \eqref{eq4}) intersecting at ridges where derivatives of $q(\lambda)$ exhibit discontinuities; in particular, $q(\lambda)$ is non-differentiable at $\lambda^{*}$. 

\textcolor{black}{For further discussion, the following definitions will be helpful:}
\textcolor{black}{\begin{definition} A vector $g(\lambda^k) \in \mathbb{R}^m$ is a subgradient of $q: \mathbb{R}^m \rightarrow \mathbb{R}$ at $\lambda^k \in \Omega$ if for all $\lambda \in \Omega$, the following holds: $q(\lambda) \leq q(\lambda^k) + \big(g(\lambda^k)\big)^{T}\cdot(\lambda - \lambda^k).$
\end{definition}}

\noindent \textcolor{black}{Note that a subgradient can exist even when $q$ is not differentiable at  $\lambda^k$. Moreover, there can be more than one subgradient of a function $q$ at a point $\lambda^k$. A way to understand the significance of the subgradient is that $q(\lambda^k) + \big(g(\lambda^k)\big)^{T}\cdot(\lambda - \lambda^k)$ is a global overestimator of $q(\lambda)$.} 
\textcolor{black}{\begin{definition} The set of subgradients of $q$ at $\lambda^k$ is referred to as the \textit{subdifferential} of $q$ at $\lambda^k$, and is denoted as $\partial q(\lambda^k)$. 
\end{definition}}

\begin{figure}[ht]
  \centering
    \includegraphics[trim=0 0 0 0,  width=0.45\linewidth, scale=0.45]{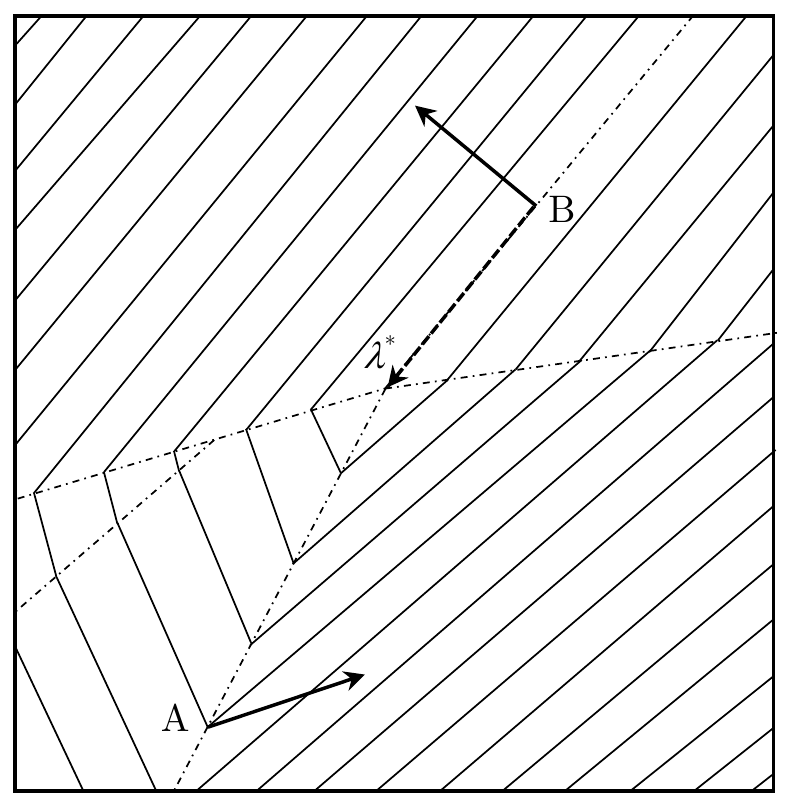}
    \caption{An example of a dual function that illustrates the difficulties associated with subgradient methods. Solid lines denote the level curves, dash-dotted lines denote the ridges of the dual function whereby the gradients are not defined (possible subgradient directions at points A and B are shown by solid arrows), and dashed lines denote the subgradient direction from point B toward optimal multipliers. This Figure is taken from \citep{bragin2022surrogate} with permission.}
    \label{fig_ex1}
\end{figure}

\noindent \textcolor{black}{While the gradient gives a direction along which a function increases the most rapidly, as discussed above, subgradients are not unique and do not necessarily provide a direction along which a non-differentiable function increases (see a subgradient direction at point A in  Figure \ref{fig_ex1}). 
Moreover, subgradients may almost be perpendicular to the directions toward optimal multipliers $\lambda^*$ thereby leading to zigzagging of $\lambda$ across ridges of the dual function (see a corresponding subgradient direction at point B in Figure \ref{fig_ex1} for illustrations).}

While Lagrangian multipliers $\lambda$ are generally fixed parameters within \eqref{eq4}, $\lambda$ are ``dual'' decision variables with respect to the dual problem \eqref{eq3}. Traditionally, \eqref{eq3} is maximized by iteratively updating  $\lambda$ by making a series of steps $s^k$ along subgradients $g(x^k,y^k)$ as:
\begin{flalign}
& \lambda^{k+1} = \lambda^k + s^k \cdot g(x^k,y^k),   \label{eq5}
\end{flalign}
where $\{x^k,y^k\}$ 
is a concise way to denote an optimal solution $\{x^*(\lambda^k),y^*(\lambda^k)\}$ to the relaxed problem \eqref{eq4} with multipliers equal to $\lambda^k.$ Within Lagrangian Relaxation, subgradients are defined as levels of constraint violations $g(x^k,y^k) \equiv \sum_{i=1}^I A_i^x \cdot x_i^k + \sum_{i=1}^I A_i^y \cdot y_i^k - b$. \textcolor{black}{Technically, since $\{x^k,y^k\}$ is a function $\lambda^k$ as discussed above, the following notation is appropriate $g(x^k,y^k) = g(\lambda^k)$. Both notations can be used interchangeably and will be denoted simply as $g^k$ for compactness as appropriate.}

\textcolor{black}{To see why $g^k$ is indeed a subgradient, consider the following: 
\begin{flalign}
& q(\lambda) \leq L(x,y,\lambda), \forall (x,y) \in \mathcal{F}.  \label{eq5a}
\end{flalign}
which is true in view of \eqref{eq4}. Adding and subtracting $\left(g^k\right)^T \cdot (\lambda^k)$ leads to
\begin{flalign}
& q(\lambda) \leq L(x,y,\lambda^k) + \big(g^k\big)^T \cdot (\lambda - \lambda^k), \forall (x,y) \in \mathcal{F}.   \label{eq5b}
\end{flalign}}
\noindent \textcolor{black}{Since the above inequality holds for all feasible $x$ and $y$, then it also holds for specific feasible values such as $x^k$ and $y^k.$ Therefore, 
\begin{flalign}
& q(\lambda) \leq L(x^k,y^k,\lambda^k) + \big(g^k\big)^T \cdot (\lambda - \lambda^k).   \label{eq5c}
\end{flalign}}
\noindent \textcolor{black}{By definition \eqref{eq4} and due to $\{x^k,y^k\} := \{x^*(\lambda^k),y^*(\lambda^k)\}$, 
\begin{flalign}
& L(x^k,y^k,\lambda^k) \equiv q(\lambda^k).   \label{eq5d}
\end{flalign}
Therefore, 
\begin{flalign}
& q(\lambda) \leq q(\lambda^k) + \big(g^k\big)^T \cdot (\lambda - \lambda^k).   \label{eq5e}
\end{flalign}
\noindent This concludes the proof.}

If inequality constraints $\sum_{i=1}^I A_i^x \cdot x_i + \sum_{i=1}^I A_i^y \cdot y_i \leq b$ are present, they are generally converted into equality constraints by introducing non-negative real-valued slack variables $z$ such that $\sum_{i=1}^I A_i^x \cdot x_i + \sum_{i=1}^I A_i^y \cdot y_i + z = b.$ Multipliers are then updated per \eqref{eq5} with subsequent projection onto the positive orthant $\{\lambda: \lambda \geq 0\}$.

Because the dual problem results from the relaxation of coupling constraints, dual values are generally less than primal values $q(\lambda^k) < f(x^k,y^k),$\footnote{In this particular case, $f(x^k,y^k) \equiv \sum_{i=1}^I \left((c_i^{x})^T \cdot x_i^k + (c_i^{y})^T \cdot y_i^k \right)$ such that $\{x^k,y^k\}$ satisfy constraints \eqref{eq2}.} i.e., there is a \textit{duality gap} - the relative difference between $q(\lambda^k)$ and $f(x^k,y^k)$.\footnote{Dual values can be used to quantify the quality of the solution $\{x^k,y^k\}$.} Because of the discrete nature of the primal problem \eqref{eq1}-\eqref{eq2}, even at optimality, the duality gap is generally non-zero, that is $q(\lambda^*) < f(x^*,y^*).$
Consequently, maximization of the dual function does not lead to an optimal primal solution $(x^*,y^*)$ or even a feasible solution. To obtain solutions feasible with respect to the original problem \eqref{eq1}-\eqref{eq2}, solutions to the relaxed problem $\{x^k_i,y^k_i\}$ are typically perturbed heuristically.\footnote{To solve MILP problems, Lagrangian Relaxation is often regarded as a heuristic. However, in dual space, the Lagrangian relaxation method is exact; the method is also capable of helping to improve solutions through the multipliers update, unlike many other heuristic methods.} Generally, the closer the multipliers are to their optimal values $\lambda^*$, the smaller the levels of constraint violations (owing to the concavity of the dual function), and, therefore, the easier the search for feasible solutions.

In short summary, the roadblocks on the way of Lagrangian Relaxation to efficiently solve MILP problems are the following: 
\begin{enumerate}
\item Non-differentiability of the dual function; 
\begin{enumerate}
\item Subgradient directions are non-ascending;
\item Necessary and sufficient conditions for extrema are inapplicable; 
\end{enumerate}
\item High computational effort is required to compute subgradient directions if the number of subsystems is large;
\item Zigzagging of multipliers across ridges of the dual function leading to many iterations required for convergence; this difficulty follows from the non-differentiability of the dual function, but this difficulty deserves a separate resolution; 
\item Solutions $\{x^k_i,y^k_i\}$ to the relaxed problem, when put together, do not satisfy constraints \eqref{eq2}. 
Moreover, 
\item Even at $\lambda^*$, levels of constraint violations may be large, and the heuristic effort to ``repair'' the relaxed problem solution $\{x^k_i,y^k_i\}$ may still be significant.    
\end{enumerate}

In order to resolve difficulty 1(b), stepsizes need to approach zero (yet, this condition alone is not sufficient), as will be discussed in the subsection that follows. This requirement puts a restriction on the methods that will be reviewed.


\textbf{Scope.} By examining the above-mentioned difficulties (which will also be referred to as D1(a), D1(b), D2, D3, D4, and D5), several stages in the development of Lagrangian Relaxation and its applications to optimizing non-smooth dual functions and solving MILP problems will be chronologically reviewed, along with specific features of the methods that address the above difficulties. 
In view of the above difficulties such as D1(b) and D2, several research directions, though having their own merit, will be excluded,
for example: 
\begin{enumerate}
\item \textbf{The Method of Multipliers.} The Alternate Direction Method of Multipliers (ADMM), which is derived from Augmented Lagrangian Relaxation (ALR) (the ``Method of Multipliers''), introduces quadratic penalties to penalize violations of relaxed constraints, improving the convergence of Lagrangian Relaxation. The two methods (ADMM and ALR), however, only converge when solving continuous primal problems. Without stepsizes approaching zero, neither method converges [in the dual space] when solving discrete primal problems and does not resolve D1(b). Nevertheless, the penalization idea underlying ALR led to the development of other LR-based methods with improved convergence as described in subsection \ref{sec27}.  

\item \textbf{The Bundle Method.} The Bundle Method's idea is to obtain the so-called $\varepsilon-$ascent direction to update multipliers \citep{zhao2002new}. Considering that the non-differentiability of dual functions may generally result in non-ascending subgradient directions, the Bundle method resolves D1(a). Since the relaxed problems need to be solved several times \citep{zhao2002new}, however, the effort required to obtain multiplier-updating directions exceeds that required in subgradient methods, thus the method does not resolve D2. 

\end{enumerate} 

\subsection{
1960's: 
Minimization of ``Unsmooth Functionals''}
\label{sec21}
Optimization of non-smooth convex functions, a direction that stems from the seminal work of \cite{polyak1967general}, is a broader subject than the optimization of $q(\lambda)$ within Lagrangian Relaxation. To present the underlying principles that support Lagrangian Relaxation to efficiently solve MILP problems, the work of \cite{polyak1967general} is discussed next. 

\noindent \textbf{Subgradient Method with ``Non-Summable'' Stepsize.} While subgradients are generally non-descending (non-ascending) for minimization (maximization) problems \cite[p.33]{polyak1967general}, convergence to the optimal solution optimizing a non-smooth function (e.g., to $\lambda^*$ maximizing $q(\lambda)$) was proven under the following (frequently dubbed as \textit{non-summable}) stepsizing formula satisfying the following conditions: 
\begin{flalign}
&  s^k > 0, \quad \lim_{k \rightarrow \infty} s^k = 0, \quad \sum_{k=1}^{\infty}{s^k} = \infty.  \label{eq6}
\end{flalign}

\noindent \textbf{Subgradient Method with Polyak Stepsize.} 
As Polyak noted in his later work \cite[p.15]{polyak1969minimization}, non-summable stepsizes lead to very slow convergence. \textcolor{black}{Intending to achieve linear\footnote{Superlinear convergence is also possible, however, 1. A reformulation of the dual problem \citep{charisopoulos2022superlinearly} is required; 2. Within the Lagrangian Relaxation framework, a dual function is generally unavailable as argued in subsection \ref{sec23}.} rate of convergence so that $\|\lambda^k - \lambda^*\|$ is monotonically decreasing, Polyak developed a stepsizing formula, which can be presented in the following way:}  
\begin{flalign}
& s^k = \gamma \cdot \frac{q(\lambda^{*}) - q(\lambda^k)}{\big\|g(x^k,y^k)\big\|^2}, 0 < \varepsilon_1 \leq \gamma \leq 2 - \varepsilon_2, \varepsilon_2 > 0.  \label{eq7}
\end{flalign}
Assuming that the function $q(\lambda)$ is strongly convex ($q(\lambda^*) - q(\lambda^k) \geq m \cdot \|\lambda^* - \lambda^k\|^2$) and satisfies the Lipshitz condition ($g^k \leq M \cdot \|\lambda^k - \lambda^*\|$) (both conditions are stated in \cite[p. 17]{polyak1969minimization}),
a rendition of Polyak's result can be presented as follows. First, consider a binomial expansion of $\|\lambda^*-\lambda^{k+1}\|^2$ as
\begin{flalign}
& \|\lambda^*-\lambda^{k+1}\|^2 = \|\lambda^*-\lambda^{k}\|^2 - 2 \cdot s^k \cdot (g^k)^T \cdot (\lambda^*-\lambda^{k}) + (s^k)^2 \cdot \|g^k\|^2. \label{eq8}
\end{flalign}
Owing to the concavity of the dual function,
\begin{flalign}
& q(\lambda^{*})-q(\lambda^k) \leq (g^k) \cdot (\lambda^* - \lambda^k). \label{eq9}
\end{flalign}
Therefore, \eqref{eq8} becomes: 
\begin{flalign}
& \|\lambda^*-\lambda^{k+1}\|^2 \leq \|\lambda^*-\lambda^{k}\|^2 - 2 \cdot s^k \cdot (q(\lambda^{*})-q(\lambda^k)) + (s^k)^2 \cdot \|g^k\|^2. \label{eq10}
\end{flalign}
Given 
the Polyak stepsize formula \eqref{eq7}, equation \eqref{eq10} becomes:
\begin{flalign}
& \|\lambda^*-\lambda^{k+1}\|^2 = \|\lambda^*-\lambda^{k}\|^2 - 2 \cdot \gamma \cdot \frac{(q(\lambda^*) - q(\lambda^k))^2}{\|g^k\|^2}+ \gamma^2 \cdot \frac{(q(\lambda^*) - q(\lambda^k))^2}{\|g^k\|^2}. \label{eq11}
\end{flalign}
The strong convexity condition $q(\lambda^*) - q(\lambda^k) \geq m \cdot \|\lambda^* - \lambda^k\|^2$ implies that $- \left(q(\lambda^*) - q(\lambda^k)\right)^2 \leq - m ^2 \cdot \|\lambda^* - \lambda^k\|^4$, therefore, given that $(2 \cdot \gamma - \gamma^2) > 0$ for $0 < \varepsilon_1 \leq \gamma \leq 2 - \varepsilon_2$ equation \eqref{eq11} becomes 
\begin{flalign}
& \|\lambda^*-\lambda^{k+1}\|^2 = \|\lambda^*-\lambda^{k}\|^2 - (2 \cdot \gamma - \gamma^2) \cdot \frac{m^2 \cdot \|\lambda^*-\lambda^k\|^4}{\|g^k\|^2}. \label{eq12}
\end{flalign}
Given that the Lipchitz condition is $g^k \leq M \cdot \|\lambda^k - \lambda^*\|$, equation \eqref{eq12} becomes
\begin{flalign}
& \|\lambda^*-\lambda^{k+1}\|^2 = \|\lambda^*-\lambda^{k}\|^2 - (2 \cdot \gamma - \gamma^2) \cdot \frac{m^2 \cdot \|\lambda^*-\lambda^k\|^4}{M^2 \cdot \|\lambda^*-\lambda^k\|^2}. \label{eq13}
\end{flalign}
Therefore
\begin{flalign}
& \|\lambda^*-\lambda^{k+1}\|^2 = \|\lambda^*-\lambda^{k}\|^2 \cdot \left(1 - (2 \cdot \gamma - \gamma^2) \cdot \frac{m^2}{M^2}\right). \label{eq14}
\end{flalign}
Within \eqref{eq7}-\eqref{eq14} and thereafter in the paper, the standard Euclidean norm will be used (unless specified otherwise). 

The Polyak stepsize \eqref{eq7} can be regarded as a creative workaround of D1(a) in the sense that a more computationally difficult problem of obtaining ascending directions at every iteration (as in the Bundle method) to ensure convergence is replaced with a provably easier problem of reducing $\|\lambda^*-\lambda^{k}\|$ at every step to guarantee convergence. 


Throughout the decades following 1969, research on the two distinct stepsize directions, i.e., ``non-summable" \eqref{eq6} and ``Polyak" \eqref{eq7}, has continued, with various groups focusing on one or the other. Both directions 
guarantee convergence to $\lambda^*$ that maximizes the dual function $q(\lambda)$ (thereby resolving D1(b)), although, up to this stage in the discussion, convergence by using Polyak stepsize is purely theoretical, since optimal dual value $q(\lambda^*)$ required within \eqref{eq7} is unknown. 

\textcolor{black}{Before we delve into the resolution of the above difficulties, it is essential to reflect upon the foundational research that has shaped our understanding of Lagrangian Relaxation for discrete programming problems. Without this understanding of seminal research that significantly influenced the trajectory of research in discrete programming, our historical exploration of the subject would remain incomplete. Even though the optimal dual value $q(\lambda^{*})$ had not been known, both the geometric convergence potential offered by the formula \eqref{eq7} and the exponential reduction of complexity upon decomposition within Lagrangian Relaxation were extensively exploited and gave rise to the resolution of MILP problems as will be discussed next.} 


\subsection{\textcolor{black}{1970's and 1980's: Application of LR and Non-Smooth Optimization to Mathematical Programming and Operations Research Problems}} \label{sec22a}

\textcolor{black}{The decades of the 1970's and 1980's were characterized by an expansive wave of both foundational research and practical applications within the domain of Lagrangian Relaxation. The fundamental research, exemplified by the works such as \cite{shapiro1971generalized}, \cite{geoffrion1974lagrangian} and \cite{fisher1974constructive}, laid the groundwork for the future development of the field. Concurrently, the field also witnessed a significant surge in the exploration and application of non-smooth optimization techniques. Notably, the Polyak method, as per equation \eqref{eq7}, emerged as a widely employed technique for solving a broad spectrum of discrete optimization problems, encompassing both pure integer and mixed-integer problems. Since the optimal dual value $q(\lambda^{*})$ within \eqref{eq7} of the associated dual value is generally unknown, to compute stepsizes, $q(\lambda^{*})$ was estimated, for example, through a ``target'' value $\overline{q}$ (e.g., \cite{held1971traveling}) or by using a feasible cost of the primal problem (e.g., \cite{fisher1976dual}). Notable applications of the method are the \textit{traveling salesman problems} (\cite{held1970traveling}, \cite{held1971traveling}), \textit{scheduling problems} (\cite{fisher1973optimal}, \cite{fisher1976dual}, \cite{muckstadt1977application}, \cite{shepardson1980lagrangean}), \textit{location problems} (\cite{cornuejols1977exceptional}, \cite{erlenkotter1978dual}) and many others. Excellent summaries of early applications of Lagrangian relaxation for discrete programming problems can be found at \cite{fisher1981lagrangian} and \cite{fisher1985applications}. Through the retrospective examination of the above research, understanding the foundational principles, insights, and techniques that have shaped the development of Lagrangian Relaxation methods in discrete programming offers us valuable perspectives as we delve deeper into the advanced methodologies and innovative applications discussed later in our survey.} 

\textcolor{black}{Transitioning into the 1990's, a transformative approach - the Subgradient-Level Method - was introduced to tackle the practical convergence issues often encountered with the Polyak method.}

\subsection{
The 1990's: 
The Subgradient-Level Method} \label{sec22} The Subgradient-Level Method \citep{goffin1999convergence} addressed difficulties associated with the unavailability of optimal [dual] value, which is needed to compute Polyak stepsize \eqref{eq7} through adaptively adjusting a level estimate based on the detection of sufficient descent and oscillations of the [dual] solutions.  

In terms of the problem \eqref{eq3}, the procedure of the method is explained as follows: the ``level'' estimate $q^{k}_{lev} = q^{k_j}_{rec} + \delta_j$ is used in place of the optimal dual value $q(\lambda^{*})$, where $q^{k}_{rec}$ is the best (highest) dual value (``record objective value'') obtained up to an iteration $k,$ and $\delta_j$ is an adjustable parameter with $j$ denoting the $j^{th}$ update of $q^{k}_{lev}.$ The main premise behind this is when $\delta_j$ is ``too large,'' then multipliers will exhibit oscillations while traveling significant (predefined) distance $R$ without improving the ``record'' value. In this case, the parameter $\delta_j$ is updated as $\delta_{j+1} = \beta \cdot \delta_j$ with $\beta = \frac{1}{2}.$ On the other hand, if $\delta_j$ is such that the dual value is sufficiently increased: $q(\lambda^k) \geq q^{k}_{lev} + \tau \cdot \delta_j,$ with $\tau = \frac{1}{2},$ then the parameter $\delta_j$ is unchanged and the distance traveled by multipliers is reset to 0 to avoid premature reduction of $\delta_j$ by $\beta$ in future iterations.  

\textcolor{black}{This method ushered in a new era of ``level-based'' optimization for non-smooth functions, which will be discussed to ultimately enhance the efficiency of the resolution of MILP problems. This advancement laid a robust foundation that has continued to resonate in subsequent methodologies and will reverberate throughout the later sections of this survey.}

Followed by an examination of resolutions of D2 and D3, further discussions of the implementation of Polyak stepsize to resolve D1 will be deferred to future subsections. 

\subsection{Fundamental Difficulties of Subgradient Methods} \label{sec23}

\noindent \textbf{High Computational Effort (D2).} In the methods reviewed thus far, non-smooth functions 
are assumed to be given. 
However, a dual function $q(\lambda)$ cannot be obtained 
computationally efficiently.
In fact, even for a given value of multipliers $\lambda^k$, minimization within \eqref{eq4} to obtain a dual value $q(\lambda^k)$ and the corresponding subgradient is time-consuming. Even then, only one possible value of the subgradient can generally be obtained; a complete description of the subgradient is generally non-attainable \citep{goffin1977convergence}. 

\noindent \textbf{Zigzagging of Multipliers (D3).} As hypothesized by \citep{goffin1977convergence}, the slow convergence of subgradient methods is due to ill-conditioning. The condition number $\mu$ is formally defined as \citep{goffin1977convergence}: $\mu = \inf\{\mu(\lambda): \lambda \in \mathbb{R}^m/P\},$ where $P = \{\lambda \in \mathbb{R}^m: q(\lambda) = q(\lambda^{*})\}$ (the set of optimal solutions) and $\mu(\lambda) = \min_u {\frac{u^T \cdot (\lambda^{*}-\lambda)}{\|u^T\| \cdot \|\lambda^{*}-\lambda\|}}$ (the cosine of the angle that subgradient form with directions toward the optimal multipliers).\footnote{Upon visual examination of the dual function illustrated in Figure \ref{fig_ex1}, the condition number is likely 0 since the subgradient emanating from point B appears to form a right angle with the direction toward optimal multipliers.} It was then confirmed experimentally when solving, for example, scheduling \cite[Fig. 3(b), p. 104]{czerwinski1994scheduling} as well as power systems problems \cite[Fig. 4, p. 774]{guan1995nonlinear} that the ill-conditioning leads to the zigzagging of multipliers across the ridges of the dual function.

To address these two difficulties, the notions of ``surrogate,'' ``interleaved'' and ``incremental'' subgradients, which do not require relaxed problems to be fully optimized to speed up convergence, emerged in the late 1990s, and early 2000s as reviewed next. 

\subsection{The late 1990's: 
Interleaved- and Surrogate-Subgradient Methods} \label{sec24} Within the Interleaved Subgradient method proposed by \citep{kaskavelis1998efficient}, multipliers are updated after solving one subproblem at a time
\begin{flalign}
& \min_{(x_i,y_i)} \left\{(c_i^x)^T \cdot x_i + (c_i^y)^T \cdot y_i + \lambda^T \cdot \big(A_i^x \cdot x_i + A_i^y \cdot y_i\big), \left\{x_i,y_i\right\} \in \mathcal{F}_i\right\},
\label{eq15}
\end{flalign}
rather than solving all the subproblems as in subgradient methods. This significantly reduces computational effort, especially for problems with a large number of subsystems. 

The more general Surrogate Subgradient method with proven convergence was then developed by \cite{zhao1999surrogate} whereby the exact optimality of the relaxed problem (or even subproblems) is not required. As long as the following ``surrogate optimality condition'' 
\begin{flalign}
& L(\tilde{x}^k,\tilde{y}^k,\lambda^k) < L(\tilde{x}^{k-1},\tilde{y}^{k-1},\lambda^k)  \label{eq16}
\end{flalign}
is satisfied, the multipliers are updated as 
\begin{flalign}
& \lambda^{k+1} = \lambda^k + s^k \cdot g(\tilde{x}^k,\tilde{y}^k),   \label{eq17}
\end{flalign}
by using the following formula
\begin{flalign}
& 0 < s^k < \gamma \cdot \frac{q(\lambda^{*}) - L(\tilde{x}^k,\tilde{y}^k,\lambda^k)}{\left\|g(\tilde{x}^k,\tilde{y}^k)\right\|^2}, \;\; \gamma < 1.  \label{eq18}
\end{flalign}
The convergence to $\lambda^{*}$ is guaranteed \citep{zhao1999surrogate}. 
Unlike that in the original Polyak formula \eqref{eq7}, parameter $\gamma$ is less than 1 to guarantee that $q(\lambda^{*}) > L(\tilde{x}^k,\tilde{y}^k,\lambda^k)$ so that the stepsize formula \eqref{eq18} is well-defined in the first place, as proven in \cite[Proposition 3.1, p. 703]{zhao1999surrogate}. Here, ``tilde'' indicates that the corresponding solutions do not need to be necessarily subproblem-optimal. Solutions $\{\tilde{x}^k,\tilde{y}^k\}$ form a set $\mathcal{S}(\tilde{x}^{k-1},\tilde{y}^{k-1},\lambda^k) \equiv \{(x,y): L(x,y,\lambda^k) < L(\tilde{x}^{k-1},\tilde{y}^{k-1},\lambda^k)\}$. Once a member $\{\tilde{x}^k,\tilde{y}^k\} \in \mathcal{S}(\tilde{x}^{k-1},\tilde{y}^{k-1},\lambda^k)$ is found, i.e., the surrogate optimality condition \eqref{eq16} is satisfied, the optimization of the relaxed problem stops and multipliers are updated per \eqref{eq17}. A case when $\mathcal{S}(\tilde{x}^{k-1},\tilde{y}^{k-1},\lambda^k) =  \emptyset$ indicates that for given $\lambda^k$ no solution better than $\{\tilde{x}^{k-1},\tilde{y}^{k-1}\}$ can be found indicating that $(\tilde{x}^{k-1},\tilde{y}^{k-1}) = (x^*(\lambda^k),y^*(\lambda^k))$ is subproblem-optimal for given value $\lambda^k$, and the multipliers are updated by using a subgradient direction.
        
The convergence proof is quite similar to that presented in \eqref{eq8}-\eqref{eq14}. The only caveat is that  $L(\tilde{x}^k,\tilde{y}^k,\lambda^k)$ is strictly speaking not a function; unlike the dual function $q(\lambda^k)$, $L(\tilde{x}^k,\tilde{y}^k,\lambda^k)$ can take multiple values for given $\lambda^k$. Therefore, the analogue of \eqref{eq9} cannot follow from concavity of $L(\tilde{x}^k,\tilde{y}^k,\lambda^k)$. It follows from the fact that the surrogate dual value is obtained without solving all the subproblems, hence 
\begin{flalign}
& q(\lambda^*) \leq L(\tilde{x}^k,\tilde{y}^k,\lambda^*). \label{eq19}
\end{flalign}     
Adding and subtracting $g(\tilde{x}^k,\tilde{y}^k)^T \cdot \lambda^k$ from the previous inequality leads to 
\begin{flalign}
& q(\lambda^*) - L(\tilde{x}^k,\tilde{y}^k,\lambda^k) \leq g(\tilde{x}^k,\tilde{y}^k)^T \cdot (\lambda^* - \lambda^k). \label{eq20}
\end{flalign}     
The procedure described in \eqref{eq10}-\eqref{eq14} follows analogously.   

In addition to the reduction of computational effort, a concomitant reduction of multiplier zigzagging has been also observed. Indeed, with an exception of the aforementioned situation whereby $\mathcal{S}(\tilde{x}^{k-1},\tilde{y}^{k-1},\lambda^k) =  \emptyset$, a solution to one subproblem \eqref{eq15} is sufficient to satisfy \eqref{eq16}. In this case, only one term within each summation of surrogate subgradient $\sum_{i=1}^I A_i^x \cdot x_i + \sum_{i=1}^I A_i^y \cdot y_i - b$ will be updated thereby preventing surrogate subgradients from changing drastically, and from zigzagging of multipliers as the result.  

\subsection{
The Early 2000's: 
Incremental Subgradient Methods} \label{sec25}

In the incremental subgradient method, a subproblem $i$ is solved before multipliers are updated (similar to the interleaved method). However, as opposed to updating all multipliers at once, the incremental subgradient method updates multipliers incrementally. After the $i^{th}$ subgradient component is calculated, multipliers are updated as 
\begin{flalign}
& \psi^{k}_i = \psi^k_{i-1} + s^k \cdot \left(A_i^x \cdot x_i^k + A_i^y \cdot y_i^k - \beta_i \right).   \label{eq21}
\end{flalign}
Here $\beta_i$ are the vectors such that $\sum_{i=1}^I \beta_i = b$, for example, $\beta_i = \frac{b}{I}.$  Only after all $i$ subproblems are solved, are the multipliers ``fully'' updated as  
\begin{flalign}
& \lambda^{k+1} = \psi^{k}_I.   \label{eq22}
\end{flalign}
Convergence results of the subgradient-level method \citep{goffin1999convergence} have been extended for the subgradient method. Variations of the method were proposed with $\beta$ and $\tau$ belonging to an interval $[0,1]$ rather than being equal to $\frac{1}{2}.$ Moreover, to improve convergence, rather than using constant $R,$ a sequence of $R_l$ was proposed such that $\sum_{l=1}^{\infty} R_l = \infty.$ 

\subsection{
2010's: 
The Surrogate Lagrangian Relaxation Method} \label{sec26}

Based on contraction mapping, Surrogate Lagrangian Relaxation (SLR) \citep{bragin2015convergence} overcomes the difficulty associated with the lack of knowledge about the optimal dual value. At consecutive iterations, the distance between multipliers must decrease, i.e., 
\begin{flalign}
& \left\| \lambda^{k+1} - \lambda^k \right\| \leq \alpha_k \cdot \left\| \lambda^{k} - \lambda^{k-1} \right\|, \quad 0 \leq \alpha_k \leq 1. \label{eq23}
\end{flalign}
Based on \eqref{eq17} and \eqref{eq23}, the stepsize formula has been derived: 
\begin{flalign}
& s^k = \alpha_k \cdot \frac{s^{k-1}\left\|g(\tilde{x}^{k-1},\tilde{y}^{k-1})\right\|}{\left\|g(\tilde{x}^k,\tilde{y}^k)\right\|}. \label{eq24}
\end{flalign}
Moreover, a specific formula to set $\alpha_k$ has been developed to guarantee convergence: 
\begin{flalign} 
& \alpha_k = 1-\frac{1}{M \cdot k^{1-\frac{1}{k^r}}}, \; M \geq 1, \; 0 
\leq r \leq 1. \label{eq25}
\end{flalign}
Since $\alpha_k \rightarrow 1,$ stepsizes within SLR are ``non-summable.'' Linear convergence can only be guaranteed outside of a neighborhood of $\lambda^*$ \cite[Proposition 2.5, p. 187]{bragin2015convergence}.

When multipliers approach their optimal values\footnote{A quality measure to quantify the quality of multipliers (i.e., how close the multipliers are to their optimal values) will be discussed in subsection \ref{sec28}}, surrogate subgradients become closer to actual subgradients, leading to reduced constraint violations due to the concavity of dual functions. This indicates that subproblems are well-coordinated, and relaxed problem solutions are near feasible ones. As a result, only a few subproblems cause infeasibility. This drastically helps with the resolution of Difficulty D4 (``solutions to the relaxed problem, when put together may not constitute a feasible solution to the original problem''): an ``automatic'' procedure to identify and ``repair'' a few subproblem solutions that cause the infeasibility of the original problem has been developed by Bragin et al. within \cite{bragin2018scalable}. 


\subsection{
The Late 2010's - Early 2020's: 
Further Methodological Advancements} \label{sec27}


\textbf{Surrogate Absolute-Value Lagrangian Relaxation \citep{bragin2018scalable}.} The Surrogate Absolute-Value Lagrangian Relaxation (SAVLR) method is designed to guarantee convergence and to speed up the reduction of constraint violations while avoiding nonlinearity and nonconvexity that would have occurred if traditional quadratic terms had been used. In the SAVLR method, the following dual problem is considered:
\begin{flalign}
& \max_{\lambda} \{q_{\rho}(\lambda): \lambda \in \Omega \subset \mathbb{R}^m\},  \label{eq26}
\end{flalign}
where
\begin{flalign}
& q_{\rho}(\lambda) = \min_{(x,y)} \Bigg\{\sum_{i=1}^I \left((c_i^x)^T \cdot x_i + (c_i^y)^T \cdot y_i\right) + \nonumber  \lambda^T \cdot \left(\sum_{i=1}^I A_i^x \cdot x_i + \sum_{i=1}^I A_i^y \cdot y_i - b\right) + \\ & \rho \cdot \Bigg\|\sum_{i=1}^I A_i^x \cdot x_i + \sum_{i=1}^I A_i^y \cdot y_i - b\Bigg\|_1, \left\{x_i,y_i\right\} \in \mathcal{F}_i, i = 1,\dots,I \Bigg\}.  \label{eq27}
\end{flalign}

\noindent The above minimization involves the exactly-linearizable piecewise linear penalties, which penalize constraint violations thereby ultimately reducing the number of subproblems that cause infeasibility mentioned in \ref{sec26} and consequently reducing the effort required by heuristics to find primal solutions. This resolves Difficulty D5.\footnote{\textcolor{black}{The parameter $\rho$ is increased as $\rho^{k+1}=\beta\cdot\rho^k,\beta>1$ until surrogate optimality condition \eqref{eq16} is no longer satisfied, which signifies that proper ``surrogate'' subgradient directions can no longer be obtained. Moreover, with heavily penalized constraint violations, the feasibility is overemphasized and subproblem solutions may get stuck at a suboptimal solution. In these situations, $\rho$ is no longer increased, rather, violations of surrogate optimality conditions prompt the reduction of penalty coefficients $\rho^{k+1}=\rho^k/\beta,\beta>1$.}}

\noindent \textbf{Distributed and Asynchronous Surrogate Lagrangian Relaxation (DA-SLR) \citep{bragin2020distributed}.} \textcolor{black}{With the emergence of technologies that enhance distributed computation, coupled with the networking infrastructure, computational tasks can be accomplished much more efficiently by using distributed computing resources than by using a single computer.} With the assumption of a single coordinator, the DA-SLR methodology has been developed to efficiently coordinate distributed subsystems in an asynchronous manner while avoiding the overhead of synchronization. Compared to the sequential Surrogate Lagrangian Relaxation \citep{bragin2015convergence}, numerical testing shows a faster convergence (12 times speed-up to achieve a gap of 0.03\% for one instance of the generalized assignment problem).

A short summary is in order here. 
While in theory, the Polyak formula offers a linear rate of convergence, the convergence rate of the Subgradient-Level Method, however, is not discussed either in the original paper by \cite{goffin1999convergence}, or in subsequent applications of the subgradient ``level-based'' ideas (e.g., by \cite{nedic2001incremental}). Likely, because of the requirements that the multipliers need to travel, either explicitly or implicitly, an infinite distance, the linear rate of convergence cannot be achieved. While SLR-based methods avoid estimating the optimal dual value, the requirement \eqref{eq25} to avoid premature termination results in the [stepsize] non-summability. 
Ideally, the goal is to avoid multiplier zigzagging, reduce the computational effort required to obtain multiplier-updating directions and achieve linear convergence.  The first step in this direction is explained in the following subsection. 

\subsection{
Early 2020's: 
Surrogate Level-Based Lagrangian Relaxation} \label{sec28}

To exploit the linear convergence potential inherent to the Polyak stepsize formula, the Surrogate ``Level-Based'' Lagrangian Relaxation (SLBLR) method has been recently developed \citep{bragin2022surrogate}. It was proven that once the following ``multiplier divergence detection'' feasibility problem
\begin{flalign} 
&\|\lambda-\lambda^{k_j+1}\| \leq \|\lambda-\lambda^{k_j}\|, \nonumber \\
&\|\lambda-\lambda^{k_j+2}\| \leq \label{eq28} 
\|\lambda-\lambda^{k_j+1}\|,   \\
&\dots \nonumber  \\ 
&\|\lambda-\lambda^{k_j+n_j}\| \leq  \nonumber 
\|\lambda-\lambda^{k_j+n_j-1}\|,
\end{flalign}
admits no feasible solution with respect to $\lambda$ (which are the decision variables in the problem above) for some $k_j$ and $n_j$, then the ``level value'' equals to
\begin{flalign} 
& \overline{q}_{j} = \max_{\kappa \in [k_j,k_j+n_j]} \overline{q}_{\kappa,j} > q(\lambda^{*}), \label{eq29}
\end{flalign}
where 
\begin{flalign} 
& \overline{q}_{\kappa,j} = \frac{1}{\gamma} \cdot s^\kappa \cdot \|g(\tilde{x}^\kappa,\tilde{y}^\kappa)\|^2 + L(\tilde{x}^\kappa,\tilde{y}^\kappa,\lambda^\kappa). \label{eq30}
\end{flalign}

In subsequent iterations, the Polyak stepsize formula is used 
\begin{flalign} 
& s^k = \zeta \cdot \gamma \cdot \frac{\overline{q}_{j} - L(\tilde{x}^k,\tilde{y}^k,\lambda^k)}{\|g(\tilde{x}^k,\tilde{y}^k)\|^2}, \zeta < 1, k = k_{j+1},\dots,k_{j+1}+n_{j+1}-1. \label{eq31}
\end{flalign}
In essence, the above formula is used until the above feasibility problem admits no solution again, at which point the level value is reset from $\overline{q}_{j}$ to $\overline{q}_{j+1}$ and the multiplier-updating process continues. It is worth noting that the optimization problem \eqref{eq3} is the maximization problem and the quality of its solutions (Lagrangian multipliers) can be quantified through the upper bound provided by $\{\overline{q}_{j}\}$.

The assumption here is that the above feasibility problem \eqref{eq28} becomes infeasible ``periodically'' and  ``infinitely often'' thereby triggering recalculations of ``level'' values $\overline{q}_{j}$, which will decrease and approach $q(\lambda^*)$ from above. The assumption is realistic since the only sure-fire way to ensure that \eqref{eq28} is always feasible is to know the optimal dual value.  

Given the above and with the addition of ``absolute-value'' penalties to facilitate the feasible solution search, the SLBLR method addresses all the difficulties D1-D5. The method inherits features from Polyak's formula \eqref{eq7}, reduction of computational effort as well as the alleviation of zigzagging from surrogate methods \citep{zhao1999surrogate, bragin2015convergence} and the acceleration of reduction of the constraint violation from \citep{bragin2018scalable}. The decisive advantage of SLBLR is provided by the practically operationalizable use of the Polyak formula \eqref{eq7} though the efficient decision-based procedure described above to determine ``level'' values without the need for estimation or heuristic adjustment of estimates of the optimal dual value. Results reported in \citep{bragin2022surrogate} indicate that the SLBLR method solves a wider range of generalized assignment problems (GAPs) to optimality as compared to other methods. With other things being equal, the ``level-based'' stepsizing of the SLBLR method \citep{bragin2022surrogate} is more advantageous as compared to the ``non-summable'' stepsizing of the SAVLR method \citep{bragin2018scalable}. Additionally, SLBLR successfully solves other problems such as stochastic job-shop scheduling and pharmaceutical scheduling outperforming commercial solvers by at least two orders of magnitude.


The SLBLR method \citep{bragin2022surrogate} is not restricted to MILP problems since linearity is not required for the above-mentioned determination of level values. The method is modular and has the potential to support plug-and-play capabilities. For example, while applications to pharmaceutical scheduling have been tested by using fixed data \citep{bragin2022surrogate}, the method is also suitable to handle urgent requests to manufacture new pharmaceutical products since such products can be introduced into the relaxed problem ``on the fly'' and the corresponding subproblems can keep being coordinated through Lagrangian multipliers without the major intervention of the scheduler and without disrupting the overall production schedule. 

\section{Conclusions and Future Directions} \label{Conclusion}

This paper intends to summarize key difficulties encountered on a path to efficiently solve MILP problems as well as to provide a brief summary of important milestones of a more than half-a-century-long research journey to address these difficulties by using Lagrangian Relaxation.
Moreover, the most recent SLBLR method is 1. general having the potential to handle general MIP problems since linearity is not required to exploit the linear rate of convergence; 2. flexible and modular having the potential to support plug-and-play capabilities with real-time response to unpredictable and/or disruptive events such as natural hazards, operational faults, and cyber-physical events, including generator outages in power systems, receiving an urgent order in manufacturing or encountering a traffic jam in transportation. With communication and distributed computing capabilities, these events can be handled by a continuous update of Lagrangian multipliers, improving system resilience; the method is thus also suitable for fast re-optimization; 3. compatible with other optimization methods such as quantum-computing as well as machine-learning methods, which can potentially be used to further improve subproblem solving 
thereby contributing to drastically reducing the CPU time supported by the fast coordination aspect of the method. 

\section*{Acknowledgements}
This work was supported in part by the U.S. National Science Foundation under Grant ECCS-1810108.

\section*{Conflict of Interest Statement}

The author declares no conflict of interest. 
\bibliographystyle{chicago}
\bibliography{IISE-Trans}
	
\end{document}